\documentclass[10pt,a4paper,leqno]{article}  

\usepackage{amsfonts,amsmath,amssymb,amsthm}
\usepackage[latin1]{inputenc}  
\usepackage[T1]{fontenc}
\usepackage[dvips]{graphicx}
\usepackage{enumerate}

\newcommand{\R}{\mathbb{R}}    
\newcommand{\C}{\mathbb{C}}

\newcommand{\bD}{\mathbb{D}}

\renewcommand{\Re}{\operatorname{Re}}
\renewcommand{\Im}{\operatorname{Im}}

\newcommand{\zbar}{{\bar{z}}}
\newcommand{\wbar}{{\bar{w}}}
\newcommand{\diam}{\operatorname{diam}}

\newcommand{\cC}{\mathcal{C}}
\newcommand{\cD}{\mathcal{D}}
\newcommand{\cE}{\mathcal{E}}
\newcommand{\cF}{\mathcal{F}}
\newcommand{\cH}{\mathcal{H}}
\newcommand{\cI}{\mathcal{I}}
\newcommand{\cJ}{\mathcal{J}}

\newcommand{\cO}{\mathcal{O}}
\newcommand{\cR}{\mathcal{R}}
\newcommand{\cS}{\mathcal{S}}

\newcommand{\D}{\partial}
\newcommand{\DD}{\nabla}

\newcommand{\loc}{\mathrm{loc}}

\let\swap=\phi
\let\phi=\varphi
\let\varphi=\swap
\let\swap=\epsilon
\let\epsilon=\varepsilon
\let\varepsilon=\swap
\let\swap=\leq
\let\leq=\leqslant
\let\leqslant=\swap
\let\swap=\geq
\let\geq=\geqslant
\let\geqslant=\swap

\theoremstyle{plain}
\newtheorem{thm}{Theorem}[section]
\newtheorem{lem}[thm]{Lemma}

\theoremstyle{definition}

\newtheorem{rem}[thm]{Remark}

\numberwithin{equation}{section}

\title{On reduced Beltrami equations and \\ linear families of quasiregular mappings}
\author{Jarmo Jääskeläinen}
\date{}
\begin{document}
\maketitle
\frenchspacing

\begin{abstract}\noindent This paper studies linear classes of planar quasiregular mappings. We give a positive answer to a conjecture of K. Astala, T. Iwaniec, and G. Martin (2009) on reduced Beltrami equations. Moreover, we use it to prove a Wronsky-type theorem for general linear Beltrami systems. This is a key to show that the associated Beltrami equation of a linear quasiregular family is unique.  
\end{abstract}

\let\thefootnote\relax\footnotetext{2010 AMS Mathematics Classification Numbers: Primary 30C62; Secondary 35J15.}

\let\thefootnote\relax

\let\thefootnote\relax\footnotetext{The author was supported by the Academy of Finland, project no. 1134757, the Finnish Centre of Excellence in Analysis and Dynamics Research, EU Research Training Network CODY Finland, and the Vilho, Yrjö and Kalle Väisälä Foundation.}

\section{Introduction}

Distinctive for the {\em reduced Beltrami equation}
\begin{equation}\label{rdc}
  \frac{\D f}{\D \zbar} = \lambda(z)\Im \biggl(\frac{\D f}{\D z} \biggr),
   \qquad |\lambda(z)| \leq k < 1, 
\end{equation} 
for almost every $z\in\Omega$, is that its solutions create an $\R$-linear family of quasiregular mappings.
The reduced equation arises naturally in a great variety of topics; for instance, in the study of linear families of quasiregular mappings, and also in the Sto\"ilow factorization and the $G$-closure problems for the general Beltrami equation
\begin{equation}\label{gen}
  \frac{\D f}{\D \zbar} = \mu(z) \frac{\D f}{\D z} + \nu(z)\overline{\frac{\D f}{\D z}\,}, \qquad |\mu(z)| + |\nu(z)|   \leq k < 1,  
\end{equation}
for almost every $z\in\Omega$. 

It is clear from the definition that a differential constraint in \eqref{rdc} is stronger than the one in the classical Beltrami equation. Hence $f$ is $K$-quasiregular with $K = \frac{1 + k}{1 - k}$. In this paper we assume an appropriate regularity, i.e., $f\in W^{1,2}_{\loc}(\Omega)$ for a domain $\Omega\subset\C$. 

Studies of the reduced Beltrami equation \eqref{rdc} indicate that for its solutions, {\em reduced quasiregular mappings}, the null Lagrangian
\[
\cJ(z,f) =  \Im \biggl( \frac{\D f}{\D z} \biggr)
\]
has many properties similar to the Jacobian determinant of a Sobolev function. We prove the following key result in this direction, answering in positive Conjecture 6.3.1 in \cite{AIM}.
\begin{thm}\label{imnon}
Suppose $f: \Omega \to \C$, $f\in W^{1,2}_{\loc}(\Omega)$, is a solution to  the reduced Beltrami equation \eqref{rdc}.
Then either $\D_z f$ is a constant or else
$$
 \Im \biggl(\frac{\D f}{\D z} \biggr) \not= 0 \qquad \text{almost everywhere in $\Omega$.}
$$
Thus if\, $\Im(\D_z f)$ vanishes on a set of positive measure, then $f(z) = az + b$, where $a\in\R$ and $b\in\C$.
\end{thm}

In Geometric Function Theory the above fact has a similar important role  as that of the Jacobian determinant of a general quasiregular map, $J(z, f) = |\D_z f|^2 - |\D_{\zbar} f|^2 \neq 0$. The null Lagrangian $\Im (\D_z f)$ appears naturally as the denominator in algebraic fractions; for example, in the $G$-compactness studies of $K$-quasiregular families.

For special properties of the reduced Beltrami equation \eqref{rdc}, we refer the reader to the recent monograph \cite{AIM}. An early application of the reduced equation can be found in \cite{Boj}. The reduced equation has generated a considerable new-found interest, see \cite{AN}, \cite{AIM}, \cite{AJ}, \cite{G-cpt}, \cite{IKO1} \cite{IKO}, \cite{KO1}, and \cite{KO2}.

First steps in giving the positive answer to the conjecture were made by F. Giannetti, T. Iwaniec, L. Kovalev, G. Moscariello, and C. Sbordone in \cite{G-cpt}. They proved the statement for global homeomorphisms, that is, homeomorphisms of the plane $\C$, when $k < \frac{1}{2}$ in \eqref{rdc}. Next, G. Alessandrini and V. Nesi showed in \cite{AN} the assertion for global homeomorphisms. A direct, and substantially simplified, proof of this result can be found in \cite{AJ} or Theorem 6.4.1 in \cite{AIM}. 

The case of plane homeomorphisms was used, for example, to study linear families of quasiconformal mappings. Moreover, by combining it with the ideas and results developed in \cite{G-cl} and \cite{G-cpt}, one can show that the family of Beltrami differential operators is $G$-compact, see Section 16.6 in \cite{AIM}. 

In spite of the close analogy with the global homeomorphic case, our proof requires more involved methods and rigorous analysis. Difference is that we consider an arbitrary domain $\Omega \subset \C$ instead of the whole plane $\C$ in the reduced Beltrami equation \eqref{rdc}; the earlier proofs use the property that a solution is a global homeomorphism $\C \to \C$. The key element of our proof is the weak reverse Hölder inequality for the solutions to so-called adjoint equations. We use it with a smoothness at a point to derive our result.

\medskip

The reduced Beltrami equation \eqref{rdc} is closely related to linear  families  of quasiregular mappings. 
The key ingredient is the following {\em Wronsky-type theorem}. It shows that the singular set of a two-dimensional linear family of quasiregular mappings has measure zero. This was conjectured for homeomorhisms in \cite{G-cl} and proven for them in \cite{AN}, \cite{AJ}. We establish a more general theorem.
\begin{thm}\label{independent}
Suppose $\Phi, \Psi \in W^{1,2}_{\loc}(\Omega)$ are solutions to \eqref{gen}. Solutions $\Psi$ and $\Phi$ are $\R$-linearly independent if and only if complex gradients $\D_z \Phi$ and $\D_z \Psi$ are pointwise independent almost everywhere, i.e.,
$$
\cJ(\Phi, \Psi) := \Im \biggl( \frac{\D \Phi}{\D z}  \overline{\, \frac{\D \Psi}{\D z} \,} \,\biggr) \neq 0 \qquad \text{almost everywhere in $\Omega$}.
$$
\end{thm}
Above $\cJ(\Phi, \Psi)$ plays the role of Wronskian. Note that, if $\Psi$ and $\Phi$ are $\R$-linearly dependent, then $\cJ(\Phi, \Psi) \equiv 0$.

Given an $\R$-linear subspace $\cF \subset W^{1,2}_{\loc}(\Omega)$, we say  that  $\cF $ is  a {\em linear family of quasiregular mappings}, if  there is  $1 \leq K< \infty$ such that for every $f \in  \cF$ the function $f$ is $K$-quasiregular in $\Omega$. The family $\cF$ is {\em generated} by the maps $f_i$, $i = \cI$,  if 
$$
\cF =\left\{ \sum_{i \in \cI} a_i\, f_i: \; a_i \in \R \right\}
$$
for some $\R$-linearly independent quasiregular mappings $ f_i: \Omega \to \C$. We always assume that the linear family is generated by countable set of functions. It quickly follows that in case of linear families that consist of quasiconformal mappings, ${\rm dim}\, \cF \leq 2$, see \cite{G-cl}. Recall that a linear family of quasiregular mappings is not always two-dimensional, e.g., $1$-quasiregular family spanned by $f_i(z) = z^i$, $i = 1, 2, 3$.

In general, quasiregularity is not preserved under linear combinations, simple example is $f(z) = k\zbar + z$, $g(z) = k\zbar - z$. However, if we have mappings that happen to be  solutions to the same general Beltrami equation \eqref{gen}, then their linear  combinations are quasiregular. Conversely, \cite{G-cl} associates to a linear two-dimensional family $\cF $ of quasiregular mappings a general Beltrami equation of the type \eqref{gen} satisfied by every $g \in \cF$, see also Remark 16.6.7 in \cite{AIM}. We show that the associated equation is unique.

\begin{thm}\label{family} For any linear family $\cF$ of quasiregular mappings, there exists a corresponding general linear Beltrami equation
$$
\frac{\D f}{\D \zbar} = \mu(z) \frac{\D f}{\D z} + \nu(z)\overline{\frac{\D f}{\D z}\,}, \qquad |\mu(z)| + |\nu(z)|   \leq k < 1, 
$$
almost everywhere in $\Omega$, satisfied by every element $f\in\cF$\!. 

Moreover, the associated equation is unique.
\end{thm}

\section{Adjoint equation}

We begin with the adjoint equation approach, similarly as in \cite{AN} and \cite{AJ}. We study the solution $f \in W^{1,2}_{\loc}(\Omega)$ to the reduced Beltrami equation
\begin{equation}\label{rdc2}
  \D_{\zbar} f(z) = \lambda(z)\Im \bigl(\D_z f(z)\bigr),
   \qquad |\lambda(z)| \leq k < 1,
\end{equation}
for almost every $z\in\Omega$.
Let us write $f(z) = u(z) + iv(z)$, where $u$ and $v$ are real-valued; similarly notate $\lambda(z) = \alpha(z) + i\beta(z)$. 
 
Using the definition of complex partial derivatives (i.e., $2\D_z = \D_x - i \D_y$ and $2\D_{\zbar} = \D_x + i \D_y$) and taking the imaginary part of the reduced equation gives
$u_y + v_x = \beta(v_x-u_y)$,  that is,
\[
u_y = \frac{\beta-1}{\beta+1} v_x.
\]
Thus 
\begin{equation}\label{uy}
 2 \Im \bigl(\D_z f(z)\bigr) = v_x - u_y =  \frac{2}{\beta+1} v_x = \frac{2}{\beta-1}u_y.
\end{equation}
Since $|\beta(z)| \leq |\lambda(z)| \leq k < 1$, the coefficients $2/(\beta(z)\pm 1)$ in \eqref{uy} are uniformly bounded from below. Hence $\Im (\D_z f)$ and  $u_{y}$ have the same zeros. 

\medskip

For the reduced Beltrami equation \eqref{rdc2}, the derivative $u_y$ is a weak solution to the adjoint equation determined by a non-divergence type operator. More precisely, consider an operator
\begin{equation}\label{Leq}
L = \sum_{i,j = 1}^{2} \sigma_{ij}(z)\frac{\D^{2}}{\D x_{i}x_{j}},
\end{equation}
where $\sigma_{ij} = \sigma_{ji}$ are measurable and the matrix
\[
\sigma(z) = 
\left[ \begin{array}{cc}
\sigma_{11}(z) & \sigma_{12}(z)  \\
\sigma_{12}(z) & \sigma_{22}(z)  
\end{array} \right]
\]
is uniformly elliptic, 
\[
\frac{1}{K}|\xi|^{2} \leq \langle\sigma(z)\xi, \xi\rangle 
= \sigma_{11}(z)\xi_{1}^{2} + 2\sigma_{12}(z)\xi_{1}\xi_{2} + \sigma_{22}(z)\xi_{2}^{2} \leqslant K|\xi|^{2}
\]
for all $\xi\in\C$ and $z \in \Omega$. Above $K$ is the {\em ellipticity constant}.
The mapping $\omega\in L^{2}_{\loc}(\Omega)$ is a weak solution to the adjoint 
equation $L^{*}(\omega) = 0$ if 
\begin{equation}\label{adjoint}
\int_{\Omega} \omega L(\phi)dm = 0, \qquad \text{for every $\phi\in C^{\infty}_{0}(\Omega)$.}
\end{equation}

To identify $u_y$ as a weak solution to an adjoint equation of the type \eqref{adjoint}, we recall that the components of solutions $f=u + iv$ to general Beltrami systems \eqref{gen} satisfy a divergence type second-order equation, see Section 16.1.5 in \cite{AIM}; note that one gets the divergence type equations for solutions to general Beltrami systems defined in any domain $\Omega$, but getting from the divergence type equations to \eqref{gen} one needs the domain to be simply connected.

In the case of the reduced Beltrami equation \eqref{rdc2}, the component $u$ satisfies 
\begin{equation}\label{div}
\operatorname{div} A\DD u = 0, \qquad
A(z) := 
\left[ \begin{array}{cc}
1 & a_{12}(z)  \\
0 & a_{22}(z)  
\end{array} \right]\!\!,
\end{equation}
where the matrix elements are
\begin{equation}\label{a1222}
a_{12} = \frac{2 \Re (\lambda)}{1- \Im(\lambda)} = \frac{2\alpha}{1-\beta}, \qquad 
a_{22} = \frac{1 + \Im(\lambda)}{1- \Im(\lambda)}  = \frac{1+\beta}{1-\beta}.
\end{equation}
Specifically \eqref{div} means that for every $\phi\in C^{\infty}_{0}(\Omega)$ 
\begin{equation}\label{id}
0 = \int_{\Omega} \DD\phi \cdot A\DD u = \int_{\Omega} \phi_{x} (u_{x} + a_{12}u_{y}) + \phi_{y}a_{22}u_{y}.
\end{equation}
But since derivatives of smooth test functions are again test functions, we can replace $\phi$ by $\phi_{y}\in C^{\infty}_{0}(\Omega)$ in \eqref{id}. Now, a straightforward calculation shows that $u_{y}$ is a weak solution to the adjoint 
equation $L^{*}(u_{y}) = 0$, where
\begin{equation}\label{L}
L = \frac{\D^{2}}{\D x^{2}} + a_{12} \frac{\D^{2}}{\D x\D y} 
+ a_{22} \frac{\D^{2}}{\D y^{2}}
\end{equation}
and $a_{12}, a_{22}$ are given by \eqref{a1222}.
We note that the original matrix $A(z)$ is not symmetric. However, the operator $L$ in \eqref{L} can be represented by the symmetric matrix 
\[
\sigma(z) = 
\left[ \begin{array}{cc}
1 & a_{12}(z)/2  \\
a_{12}(z)/2 & a_{22}(z)  
\end{array} \right]
\]
and, as $|\lambda(z)| \leq k < 1$, from \eqref{a1222} we see that $\sigma$ is uniformly elliptic.

\section{Weak reverse Hölder inequality}

\begin{thm}\label{reverse}
Let $\omega\in L^2_{\loc}(\Omega)$ be a real-valued weak solution to the adjoint equation $L^*(\omega) = 0$ of the type \eqref{adjoint}. Then a weak reverse Hölder inequality holds for $\omega$\!; namely,
\begin{equation}\label{wrH}
\biggl(\frac{1}{r^2}\int_{B} \omega^2 dm\biggl)^{1/2} \leq \frac{c}{r^2}\int_{2B} |\omega|dm,
\end{equation}
for every disk $B := \bD(a, r)$ such that $2B := \bD(a, 2r)\subset\Omega$. The constant $c$ depends only on the ellipticity constant $K$\!.
\end{thm}
There is a stronger result for non-negative solutions: a reverse Hölder inequality holds, see \cite{FS} or Theorem 6.4.2 in \cite{AIM}.

We start with a well-known interpolation inequality, for instance, \cite{BL}, \cite{Mazja}, and \cite{Triebel}.
\begin{lem}
Let $g\in W^{2,p}(U)$ be real-valued, where $U\subset\C$ is a bounded smooth domain and $p > 1$. Then
\begin{equation}\label{inter}
\|\DD g\|_{L^{p}(U)} \leq c\bigl(\|D^2 g\|_{L^{p}(U)} + \|g \|_{L^{p}(U)}\bigr),
\end{equation}
where $c$ depends only on $U$\!.
\end{lem}

\begin{proof}
By a basic interpolation between Sobolev spaces, see \S\,4.3.1, Theorem 1 with \S\,2.4.2, (11) in \cite{Triebel},
\[
\bigl(W^{2, p}(U), W^{0, p}(U) \bigr)_{1/2} = W^{1,p}(U),
\]
and thus
\[
\| g\|_{W^{1,p}(U)} \leq c(U)\|g\|_{W^{2, p}(U)}^{1/2}\|g\|_{W^{0, p}(U)}^{1/2}.
\]
Then we use the inequality 
\[
ab \leq \frac{\epsilon}{2} a^2 + \frac{1}{2\epsilon} b^2, \qquad a,b \geq 0, \quad \epsilon > 0,
\]
to derive 
\begin{align*}
\|\nabla g\|_{L^{p}(U)} &\leq\| g\|_{W^{1,p}(U)} \\
&\leq \epsilon c(U)\bigl(\|g\|_{L^{p}(U)} + \|\nabla g\|_{L^{p}(U)} +  \|D^2 g\|_{L^{p}(U)}\bigr) + c(U, \epsilon)\|g\|_{L^{p}(U)}.
\end{align*}
Choosing $\epsilon$ sufficiently small the term with $\nabla g$ is absorbed to the left-hand side, and we arrive at the estimate \eqref{inter}.
\end{proof}

Next, let us recall some of the key estimates in this connection. There is always a unique solution to the following Dirichlet problem in a bounded domain $\cD$ with a thick boundary, see Chapter 17 in \cite{AIM},
\begin{equation}\label{dirichlet}
L(g) = h, \qquad h\in L^2(\cD), \qquad g\in W^{2,2}(\cD) \quad \text{with $g=0$ on $\D\cD$},
\end{equation}
where the operator $L$ is of the type \eqref{Leq}.
By the Alexandrov-Bakel$'$man-Pucci maximum principle, see Theorem 17.3.1 in \cite{AIM} or Theorem 9.1 in \cite{GT},
\begin{equation}\label{abp}
\|g\|_{L^\infty(\cD)} \leq c\diam(\cD)\|h\|_{L^2(\cD)}
\end{equation}
and $c$ depends only on the ellipticity constant $K$\!.
\medskip

\begin{lem}\label{gradbound}
Let $g$ be a real-valued solution to the Dirichlet problem \eqref{dirichlet}. If $L(g)= 0$ in a subdomain $V\subset\cD$\!, then for every relatively compact smooth subdomain $V' \subset V$ 
\begin{equation*}
\|\DD g\|_{L^\infty(V')} \leq c\|g\|_{L^\infty(\cD)},
\end{equation*}
where $c$ depends on $p \in \left(2, \frac{2K}{K-1}\right)$, $V$\!, and $V'$\!.
\end{lem}

\begin{proof}
We use the interior regularity from \cite{Pucci}; alternatively, see Chapter 17 in \cite{AIM}. By Lemma 4.1 in \cite{Pucci}, the complex gradient $g_z$ is quasiregular in $V$\!. Further, by Corollary 5.1 in \cite{Pucci},
\[
D^2g \in L_{\loc}^p(V), \qquad  2 \leq p < \frac{2K}{K-1}.
\]
Moreover, the corollary implies for every relatively compact subdomain $V' \subset V$\! the uniform estimate 
\begin{equation}\label{interior}
\|D^2g \|_{L^p(V')} \leq c\|g \|_{L^1(V)},
\end{equation}
where $c$ depends on $p$, $V$\!, and $V'$\!. Actually, for the potential function, $g \in W^{2,p}_{\loc}(V)$ for all $2 < p < \frac{2K}{K-1}$. Thus by the Sobolev embedding, or more strictly by Morrey's inequality, we achieve the following estimate for the Hölder norm
\[
\|\DD g \|_{C^{0, \gamma}(\overline{V'})} \leq c(p, V')\|\DD g\|_{W^{1, p}(V')} \leq  c(p, V') \bigl(\|\DD g\|_{L^{p}(V')} + \|D^2 g\|_{L^{p}(V')}\bigr),
\]
where $\gamma = 1 - 2/p$.
Now we use the interpolation inequality \eqref{inter} to write
\[
\|\DD g\|_{L^{p}(V')} \leq c(V')\bigl(\|D^2 g\|_{L^{p}(V')} + \|g \|_{L^{p}(V')}\bigr).
\]
Hence, by combining the previous estimates with \eqref{interior}, we have shown
\[
\|\DD g\|_{L^\infty(V')} \leq \|\DD g \|_{C^{0, \gamma}(\overline{V'})} \leq c(p, V, V')\|g\|_{L^\infty(\cD)}.
\]
\end{proof}

In the following proof, we will also use the consequence of the Leibniz rule,
\begin{equation}\label{leibniz}
L(\phi g) = 2\langle\sigma\DD\phi, \DD g \rangle + \phi L(g) + gL(\phi), \qquad \phi, g\in W^{2,2}(\cD),
\end{equation}
where $L$ is an operator defined as in \eqref{Leq}.

\bigskip

\begin{proof}[Proof of Theorem \ref{reverse}]
Without loss of generality we can assume $a = 0$. We show the claim for the unit disk $\bD$ and then use a rescaling argument.

For the unit disk case, it is enough to prove
\begin{equation}\label{enough}
\int_{\bD} \omega^2dm \leq c(K)\|\omega \|_{L^2(\bD)} \int_{2\bD}|\omega|dm.
\end{equation}

We solve the Dirichlet problem \eqref{dirichlet} for $\cD = 2\bD$ and $h = \omega\chi_{\bD}\in L^2(2\bD)$. As in \eqref{dirichlet} we notate the $W^{2,2}$-solution by $g$.

Let $1 < \delta <  4/3$ and $\phi\in C^{\infty}_{0}\bigl((3/2) \delta\bD\bigr)$ satisfy $\phi \equiv 1$ on $\delta\bD$ with $\bigl|\D^{\alpha} \phi / \D x_{\alpha}\bigr| \leq c_{\alpha}$ for $|\alpha|\leq 2$. Now $\phi g \in W^{2,2}_0(2\bD)$. Since $\omega$ is an adjoint solution in $2\bD$, $L^*(\omega) = 0$, by approximating with smooth functions we find that
\[
\int_{2\bD} \omega L(\phi g)dm =0.
\]

The consequence of the Leibniz rule \eqref{leibniz} gives
\begin{align*}
\int_{\bD} \omega^2 &= \int_{2\bD} \omega L(g)\phi = -2\int_{2\bD} \omega\langle\sigma\DD\phi, \DD g\rangle  -\int_{2\bD} \omega g L(\phi) \\
&\leq 2\int_{2\bD} |\omega||\langle\sigma\DD\phi, \DD g\rangle| + \int_{2\bD} |\omega| |g| |L(\phi)|\\
&\leq c(K)\int_{(3/2)\delta\bD\setminus \delta\overline{\bD}} |\omega||\DD \phi||\DD g| + \|g\|_ {L^\infty(2\bD)}\int_{2\bD}|\omega||L(\phi)|\\
&\leq c(K)\int_{(3/2)\delta\bD\setminus \delta\overline{\bD}} |\omega||\DD g| + c(K)\|g\|_ {L^\infty(2\bD)}\int_{2\bD}|\omega|.
\end{align*}

Further, by using Lemma \ref{gradbound} to the sets 
$$
V\! := 2\bD\setminus\overline{\bD} \subset 2\bD =: \cD \qquad  \text{and} \qquad  V'\! := (3/2)\delta\bD\setminus \delta\overline{\bD},
$$
we have
\begin{equation*}
\int_{\bD} \omega^2 \leq c(K)\|g\|_ {L^\infty(2\bD)}\int_{2\bD}|\omega|.
\end{equation*}
The inequality \eqref{enough} follows by the consequence of the Alexandrov-Bakel$'$man-Pucci maximum principle \eqref{abp}.

\smallskip

We are left with the rescaling. Assume $2B\subset\Omega$. We set 
$$
\omega_{r}(z) = \omega(rz) \qquad \text{and} \qquad L_r(\phi)(z) = L(\phi)(rz).
$$
By definition, $L_r$ is a uniformly elliptic non-divergence type operator with the same ellipticity constant $K$ as the operator $L$. Further, $\omega_r$ is an adjoint solution in $\Omega_r := \{z \in \Omega : rz\in\Omega \}$ for $L_r$ and $2\bD \subset \Omega_r$. Thus the above proof shows
\begin{equation}\label{dr}
\biggl(\int_{\bD} \omega_{r}^2(z) dm(z)\biggl)^{1/2} \leq c(K)\int_{2\bD} |\omega_r(z)|dm(z).
\end{equation}
Next we use the change of variables to get
\begin{equation*}
r^2\int_{\bD} \omega_{r}^{2}(z) dm(z) = \int_{B} \omega^2(z) dm(z), \qquad r^2\int_{2\bD} |\omega_r(z)|dm(z) = \int_{2B} |\omega(z)|dm(z).
\end{equation*}
Combining the above calculation with \eqref{dr} gives
\[
\biggl(\frac{1}{r^2}\int_{B} \omega^2(z) dm(z)\biggl)^{1/2} \leq \frac{c(K)}{r^2}\int_{2B} |\omega(z)|dm(z).
\]
\end{proof}

\medskip

\begin{rem}\label{impro}
It is well-known, see Gehring's lemma, for example, from Section 4.3 in \cite{BI} or Chapter 14 in \cite{gft}, that a weak reverse Hölder inequality improves integrability: if a weak reverse Hölder inequality \eqref{wrH} holds for $\omega$, then there is $p > 2$ such that 
$$
\biggl(\frac{1}{r^2}\int_{B} |\omega|^p dm\biggl)^{1/p} \leq \frac{c(K)}{r^2}\int_{2B} |\omega|dm.
$$
\end{rem}

\section{Zeros of infinite order}

A weak reverse Hölder inequality implies that almost every zero is of infinite order; that is,
\begin{thm}\label{zeros}
Let $\omega$ satisfy a weak reverse Hölder inequality \eqref{wrH}. Then, for almost every zero $z_0$ of\, $\omega$ and for every positive integer $N$\!, there is\, $r_0(z_0, N) > 0$ such that
$$
\int_{\bD(z_0, r)}|\omega|dm \leq \frac{r^N}{r_{0}^{N}}\int_{\bD(z_0, 2r_0)}|\omega|dm = \cO(r^N), \qquad 0 < r \leq r_0(z_0, N).
$$
\end{thm}

\begin{proof}
We use the iteration argument from pp. 299--300 in \cite{BI} and Theorem 14.5.1 in \cite{gft}.

Set $E = \{z\in\Omega : \omega(z) = 0\}$. Assume $|E| > 0$.
Let $z_0$ be a point of density of $E$. We fix a positive integer $N$\!. Since $z_0$ is a density point, we find that, for $r_0 := r_0(z_0, N)$ sufficiently small,
$$
|\bD(z_0, \delta r_0)\setminus E| \leq \frac{(\delta r_0)^2}{c^2\, 2^{2N}}
$$
holds for all $0 < \delta \leq 1$, where $c$ is the constant from the weak reverse Hölder inequality. Thus
\begin{align*}
\int_{\bD(z_0, \delta r_0)}|\omega| &= \int_{\bD(z_0, \delta r_0)\setminus E}|\omega| \leq |\bD(z_0, \delta r_0)\setminus E|^{1/2}\biggl(\int_{\bD(z_0, \delta r_0)}|\omega|^2\biggr)^{1/2} \\
&\leq |\bD(z_0, \delta r_0)\setminus E|^{1/2}\frac{c}{\delta r_0}\int_{\bD(z_0, 2\delta r_0)} |\omega| \\
&\leq \frac{1}{2^N}\int_{\bD(z_0, 2\delta r_0)} |\omega|.
\end{align*}

Iterating yields for $k = 1, 2, \ldots$
$$
\int_{\bD(z_0, 2^{-k}r_0)}|\omega| \leq \frac{1}{2^{(k+1)N}}\int_{\bD(z_0, 2r_0)} |\omega|.
$$
For each $0 < r \leq r_0$ there exists $k$ such that 
$$
2^{-k}r_0 \leq r < 2^{-k +1}r_0.
$$ 
Hence
\begin{align*}
\int_{\bD(z_0, r)}|\omega|&\leq \int_{\bD(z_0, 2^{-k +1}r_0)}|\omega|\leq \frac{1}{2^{kN}}\int_{\bD(z_0, 2r_0)} |\omega| \\
&= \frac{r_0^N}{2^{kN}r_{0}^{N}}\int_{\bD(z_0, 2r_0)} |\omega| \\
&\leq \frac{r^N}{r_{0}^{N}}\int_{\bD(z_0, 2r_0)} |\omega| = \cO(r^N).
\end{align*}
\end{proof}

\section{Proof of Theorem \ref{imnon}}

We are ready to answer in positive Conjecture 6.3.1 in \cite{AIM}.

\begin{proof}[Proof of Theorem \ref{imnon}] Let $f = u + iv$ be a solution to the reduced Beltrami equation
$$
  \D_{\zbar} f(z) = \lambda(z)\Im \bigl(\D_z f(z)\bigr),
   \qquad |\lambda(z)| \leq k < 1, 
$$
for almost every $z\in\Omega$. By \eqref{uy}
$$
E := \{z\in\Omega : \Im\bigl(\D_{z}f(z)\bigr) = 0\} = \{z\in\Omega : u_y(z) = 0\}.
$$ 

Assume $|E| > 0$. We have shown above that $u_y$ is a solution to the adjoint equation and thus a weak reverse Hölder inequality holds for $u_y$. Further, by Theorem \ref{zeros},
\begin{align}\label{zbarinf}
\int_{\bD(z_0, r)}|\D_{\zbar}f| &\leq k\int_{\bD(z_0, r)}|\Im\bigl(\D_{z}f\bigr)| \leq \frac{k}{1-k}\int_{\bD(z_0, r)}|u_y| = \cO(r^N),
\end{align}
for almost every $z_0\in E$ and for all positive integers $N$\!, when $r > 0$ is small enough. In the second inequality we use \eqref{uy} again. 

\medskip

\noindent {\em A. \; Series representation}\nopagebreak

\smallskip

\noindent We will prove that for almost every $z_0\in E$ and for all positive integers $n$,
\begin{equation}\label{inforder}
f(w) = c_0 + c_1\,(w - z_0) + \cE(w) \qquad \text{near the point $z_0$,}
\end{equation}
where  $c_0 \in \C$, $c_1\in\R$ are constants depending only on $f$ and $z_0$ and 
\begin{equation}\label{Derror}
\int_{\bD(z_0, r)} |D\cE| dm = \cO(r^{n+1})
\end{equation}
holds for small enough $r > 0$. We deduce the statement of our theorem from this by quasiregularity.

\medskip

Fix a positive integer $n$. Choose $z_0\in E$ and $r_0 \in (0, 1]$ such that $\bD(z_0, 2r_0)\subset\Omega$ and \eqref{zbarinf} holds for $N = n + 2$ and $0 < r \leq r_0$. By the weak reverse Hölder inequality \eqref{wrH} and  the improved integrability, see Remark \ref{impro}, \eqref{zbarinf} implies for some $p_0 > 2$
\begin{equation}\label{int2}
\int_{\bD(z_0, r)}|\D_{\zbar}f|^p dm = \cO(r^{Np - 2p + 2}), \qquad \text{$2 \leq p \leq p_0$,}
\end{equation}
when $0 < r \leq r_0$. Alternatively one could use Astala's higher integrability, see \cite{ast} or Theorem 13.2.3 in \cite{AIM}.

Suppose $w\in\bD(z_0, r_0)$. We begin by showing that 
\begin{equation}\label{inforder1}
f(w) = \sum_{j= 0}^{n-1} c_j\,(w - z_0)^j + \cE(w), \qquad \int_{\bD(z_0, r)} |D\cE| dm = \cO(r^{n+1}),
\end{equation}
where $0 < r \leq r_0$ and $c_j \in \C$ are constants depending only on $f$ and $z_0$. 

Smoothness at a point has been studied, for example, in \cite{dynkin} and we use a few similar ideas. 

\medskip

\noindent {\em Step 1. \; Generalized Cauchy formula}\nopagebreak

\smallskip

\noindent First, the generalized Cauchy formula gives
\begin{equation*}
f(w) = \frac{1}{2\pi i}\int_{\D \bD(z_0, r_0)}\frac{f(z)}{z - w}dz + \frac{1}{\pi}\int_{\bD(z_0, r_0)}\frac{\D_{\zbar}f(z)}{w - z}dm(z).
\end{equation*}
Since the first term is analytic in the disk $\bD(z_0, r_0)$, using the Taylor expansion about $z_0$ it can be written in the form 
$$
\sum_{j = 0}^{n-1} a_j\, (w - z_0)^j + R_{n}(w), \qquad R_{n}(w) =  \cO(|w - z_0|^{n}).
$$
For the second term,
\begin{align*}
\frac{1}{\pi}\int_{\bD(z_0, r_0)}\frac{\D_{\zbar}f(z)}{w - z}dm(z) &= -\sum_{j = 0}^{n-1} (w - z_0)^j\,\frac{1}{\pi}\int_{\bD(z_0, r_0)}\frac{\D_{\zbar}f(z)}{(z - z_0)^{j+1}}dm(z)\\
&\quad+  (w - z_0)^{n}\,\frac{1}{\pi}\int_{\bD(z_0, r_0)}\frac{\D_{\zbar}f(z)}{(z - z_0)^{n}(w- z)}dm(z)\\
&=: \sum_{j = 0}^{n-1} b_j\,(w - z_0)^j  + T(w), 
\end{align*}
as soon as we show the convergence of the coefficient integrals
\begin{equation}\label{coef}
|b_j| \leq \frac{1}{\pi}\int_{\bD(z_0, r_0)}\frac{|\D_{\zbar}f(z)|}{|z - z_0|^{j+1}}dm(z), \qquad j = 1, 2, \ldots, n-1.
\end{equation}
Observe that after we have the convergence of integrals, $f(w)$ is a sum of a holomorphic part and $T$\!.

\medskip

\noindent {\em Step 2. \; Convergence of the integrals}\nopagebreak

\smallskip

\noindent Dividing in annuli 
$$
A_k := \bD(z_0, 2^{-k + 1}r_0) \setminus \bD(z_0, 2^{-k}r_0),
$$
for $j = 1, 2, \ldots n-1$,
\begin{align*}
\int_{\bD(z_0, r_0)}\frac{|\D_{\zbar}f(z)|}{|z - z_0|^{j+1}}dm(z) &= \sum_{k = 1}^{\infty}\int_{A_k}\frac{|\D_{\zbar}f(z)|}{|z - z_0|^{j+1}}dm(z) \\
&\leq \sum_{k = 1}^{\infty}\frac{1}{(2^{-k}r_0)^{j+1}}\int_{\bD(z_0, 2^{-k + 1}r_0)}|\D_{\zbar}f|dm \\
&\leq c2^{n+2}\sum_{k = 1}^{\infty}\frac{1}{2^{2k}},
\end{align*}
by \eqref{zbarinf}. Indeed, our choice $N = n + 2$ gives for $k = 1, 2, \ldots$, 
$$
\frac{r_0^N}{r_0^{j+1}}\,\frac{2^{(-k + 1)N}}{2^{-k(j+1)}} \leq \frac{2^{n + 2}}{2^{2k}}, \qquad j = 1, 2, \ldots n-1.
$$
Thus the coefficient integrals \eqref{coef} converge. For a future reference, the same reasoning shows
\begin{equation}\label{fut}
\int_{\bD(z_0, r_0)}\frac{|\D_{\zbar}f(z)|^p}{|z - z_0|^{np}}dm(z) \leq c(K, p_0, n)\sum_{k = 1}^{\infty}\frac{1}{2^{2k}}, \qquad 2 \leq p \leq p_0,
\end{equation}
by \eqref{int2}.

\medskip

\noindent {\em Step 3. \; Remainder term}\nopagebreak

\smallskip

\noindent Set $c_j = a_j - b_j$ for $j = 0, \ldots, n - 1$. Thus, as the remainder term in \eqref{inforder1}, we have $\cE = R_n + T$, where $R_n$ is holomorphic with $R_n(w) = \cO(|z - w|^n)$ and 
$$
T(w) = (w - z_0)^{n}\,\frac{1}{\pi}\int_{\bD(z_0, r_0)}\frac{\D_{\zbar}f(z)}{(z - z_0)^{n}(w- z)}dm(z).
$$
To prove \eqref{inforder1}, we are left to show the estimate for the derivative $D\cE$. 

Recall $|D\cE| = |\D_{\zbar}\cE| + |\D_z\cE|$. By definition, the derivative of the holomorphic part $DR_n$ has the correct convergence rate. After combining $\D_{\zbar}f = \D_{\zbar}\cE = \D_{\zbar}T$ with \eqref{zbarinf}, we see that only the estimation of $\D_z T$ remains.

The integral term in $T$ is the Cauchy transform $\cC$ of 
$$
F(z) := \frac{\chi_{\bD(z_0, r_0)}(z)\D_{\zbar}f(z)}{(z - z_0)^{n}}, \qquad \cC(F)(w) := \frac{1}{\pi}\int_{\C}\frac{\chi_{\bD(z_0, r_0)}(z)\D_{\zbar}f(z)}{(z - z_0)^{n}(w- z)}dm(z).
$$
The above is well-defined, since $F\in L^2(\C)$ by \eqref{fut}. Now, for almost every $w$,
\begin{equation}\label{appT}
\D_z T(w) = n(w - z_0)^{n-1}\cC(F)(w) +  (w - z_0)^{n}\cS(F)(w),
\end{equation}
where the Beurling transform $\cS$ is given by the principal value integral
$$
\cS(F)(w) := - \frac{1}{\pi}\int_{\C}\frac{\chi_{\bD(z_0, r_0)}(z)\D_{\zbar}f(z)}{(z - z_0)^{n}(w - z)^2}dm(z).
$$ 

We start with the first term in \eqref{appT}. By inequality \eqref{fut}, there is a Hölder conjugate pair $1 < q < 2 < p < \infty$ such that $F\in L^p(\C) \cap L^q(\C)$, since $F$ has a compact support. Thus $\cC(F) \in C_0(\hat{\C})$, and moreover, 
$$
\| \cC(F)\|_{L^{\infty}(\C)} \leq \frac{1}{\sqrt{2 - q}}\bigl(\|F \|_{L^p(\C)} + \|F \|_{L^q(\C)}\bigr),
$$
see, for example, Theorem 4.3.11 in \cite{AIM}. We have
\begin{align*}
\int_{\bD(z_0, r)}|n(w - z_0)^{n-1}\cC(F)(w)|dm(w) &\leq cr^{n + 1}.
\end{align*}

For the second term in \eqref{appT} the Hölder inequality implies
\begin{align*}
\int_{\bD(z_0, r)}|w - z_0|^{n}|\cS(F)(w)|dm(w) &\leq \biggl(\int_{\bD(z_0, r)}|w - z_0|^{2n}dm(w)\biggr)^{1/2}\|\cS(F)\|_{L^2(\C)} \\
&\leq cr^{n + 1}\|F\|_{L^2(\C)}.
\end{align*}

We have proven the estimate \eqref{inforder1}.

\medskip

\noindent {\em Step 4. \; No higher-order terms}\nopagebreak

\smallskip

\noindent Observe for $j = 2, \ldots, n - 1$
\begin{align*}
\int_{\bD(z_0, r)} &\bigl|\Im\bigl(j\, c_j(w - z_0)^{j - 1}\bigr)\bigr|dm(w) \\ &= \frac{j\,r^{j +1}}{j + 1}\int_{0}^{2\pi}\bigl|\beta_j \cos\bigl((j - 1)\theta\bigr) + \alpha_j\sin\bigl((j - 1)\theta\bigr)\bigr|d\theta = cr^{j + 1},
\end{align*}
where we notate $c_j = \alpha_j + i\beta_j$. Further, straight from \eqref{inforder1}
$$
\Im\bigl(\D_z f(w)\bigr) = \Im c_1 + \sum_{j = 2}^{n - 1} \Im\bigl(j\,c_j(w - z_0)^{j - 1} \bigr) + \Im\bigl( \D_z \cE(w) \bigr).
$$
Now, the estimates for the convergence rate in \eqref{zbarinf} and \eqref{inforder1} imply $\Im c_1 = 0$ and $c_j = 0$, $j = 2, \ldots, n -1$.

We have shown the series representation \eqref{inforder} with the estimate \eqref{Derror}.

\medskip

\noindent {\em B. \; Conclusion by quasiregularity}\nopagebreak

\smallskip

\noindent We use similar methods as  in pp. 299--300 \cite{BI} and Section 16.10 in \cite{gft}. In the references these ideas are used to prove that the Jacobian of a nonconstant quasiregular mapping is nonvanishing almost everywhere.

The constant $c_1$ in \eqref{inforder} is real and hence $g(w) := f(w) - c_0 - c_1(w - z_0)$ solves the same reduced Beltrami equation as $f$\!. Therefore, $g$ is quasiregular with the following property for every positive integer $n$
$$
\int_{\bD(z_0, r)}|Dg|dm = \int_{\bD(z_0, r)}|D\cE|dm = \cO(r^{n+1}), \qquad 0 < r \leq r_0.
$$
Since the weak reduced Hölder inequality holds for $|Dg|$ by quasiregularity, we achieve for every positive integer $N$
$$
\biggl(\int_{\bD(z_0, r)}|Dg|^2dm\biggr)^{1/2}  \leq \cO(r^{N+1}), \qquad \text{when $r$ is small enough}.
$$

Quasiregularity and a version of Morrey's inequality implies the Hölder continuity of the form
$$
|g(z_0) - g(w)| \leq c\biggl(\frac{|z_0 - w|}{r}\biggr)^{\alpha(K)}\biggl(\int_{\bD(z_0, r)}|Dg|^2dm\biggr)^{1/2},  \qquad w\in \bD(z_0, r/2),
$$
and $0 < \alpha(K) < 1$, see, for instance, Theorem 5.2 in \cite{BI}.
Thus
\begin{equation}\label{prop}
\sup_{|z_0 - w| < \,r/2}|g(z_0) - g(w)| = \cO(r^{N + 1}).
\end{equation}
This proves our statement: $g$ is quasiregular and hence the classical Sto\"ilow factorization holds; that is, $g = h \circ G$, where $h$ is holomorphic and $G$ a quasiconformal homeomorphism. If $g$ is nonconstant, the quasisymmetry of $G$ and $h(z) = \cO(|z - G(z_0)|^m)$, $m \geq 1$, imply  that there exists $\gamma > 0$ such that
$$
cr^{\gamma} \leq \sup_{|z_0 - w| < \, r/2}|g(z_0) - g(w)|.
$$
This would be a contradiction with \eqref{prop}. Thus $g$ is a constant; and $f(z) = c_0 + c_1 (z_0 - w)$, where $c_0\in\C$, $c_1\in\R$, proving our claim.
\end{proof}

\section{Linear families of quasiregular mappings. Proofs of Theorems \ref{family} and \ref{independent}}

\begin{proof}[Proof of Theorem \ref{independent}] Suppose $\Phi, \Psi \in W^{1,2}_{\loc}(\Omega)$ are solutions to the general Beltrami equation \eqref{gen}. Moreover, assume $\Phi$, $\Psi$ are not affine combinations of each other. We show that $\D_z\Phi$ and $\D_z\Psi$ are linearly independent over the field $\R$, that is,
$$
\Im \biggl( \frac{\D \Phi}{\D z}  \overline{\, \frac{\D \Psi}{\D z} \,} \,\biggr) \neq 0 \qquad \text{almost everywhere in $\Omega$.}
$$

We can assume $\Phi$ is nonconstant. As a nonconstant quasiregular mapping, it follows that $\Phi$ is discrete, open, and the branch set consists of isolated points. Thus it is enough to study points outside the branch set. Let $z_0$ be such a point. There exists a ball $B := \bD(z_0, r)$ such that $\Phi|_B : B \to \Phi(B)$ is a homeomorphism, hence quasiconformal. From the Sto\"ilow factorization of general Beltrami equations, Theorem 6.1.1 in \cite{AIM}, we know that
$$
\Psi = F \circ \Phi \qquad \text{in $B$},
$$ 
where $F$ solves the reduced Beltrami equation \eqref{rdc} in $\Phi(B)$ with
$$
\lambda(w) = \frac{-2i\, \nu(z)}{1 + |\nu(z)|^2 - |\mu(z)|^2}, \qquad w = \Phi(z), \quad z \in B.
$$

Let $z\in B$. Using the chain rule and identities 
$$
J(z, f) h_{\wbar}(w) = - f_{\zbar}(z), \qquad J(z, f) h_{w}(w) = \overline{f_{z}(z)},
$$
where $h = f^{-1}$ and $w = f(z)$, we arrive at
\begin{align*}
J(z, \Phi) F_w(w) &= \Psi_z(z)\overline{\Phi_z(z)} - \Psi_{\zbar}(z)\overline{\Phi_{\zbar}(z)} \\
&= (1 - |\mu|^2)\Psi_z\overline{\Phi_z} - |\nu|^2\overline{\Psi_z}\Phi_z - 2\Re (\mu\overline{\nu}\Psi_z\Phi_z), \qquad w = \Phi(z).
\end{align*}
Thus
$$
J(z, \Phi)\Im(F_w \circ \Phi) = (-1 + |\mu|^2 - |\nu|^2)\Im(\Phi_z\overline{\Psi_z}).
$$
Since $\Phi|_B$ preserves sets of zero measure, the statement follows by Theorem \ref{imnon}.
\end{proof}

\bigskip

\begin{proof}[Proof of Theorem \ref{family}] \; {\em Step 1. Two-dimensional family}.\; It is known that for any linear two-dimensional family $\cF$ of quasiregular mappings $\Omega \to \C$ there exists a corresponding general Beltrami equation
\begin{equation}\label{geneq}
\frac{\D f}{\D \zbar} = \mu(z) \frac{\D f}{\D z} + \nu(z)\overline{\frac{\D f}{\D z}\,}, \qquad |\mu(z)| + |\nu(z)|   \leq k < 1, 
\end{equation}
almost everywhere in $\Omega$, satisfied by every element $f\in\cF$, see the beginning of Section 5.3 in \cite{G-cl} or the proof of Theorem 16.6.6 in \cite{AIM}. We recall the ideas of the proof for the reader's convenience, and further, show that the associated equation is unique. 

\medskip

Assume $\Phi, \Psi \in W^{1,2}_{\loc}(\Omega)$ generate a linear family $\cF$ of $K$-quasiregular mappings. The goal is to find coefficients $\mu$ and $\nu$ such that
\begin{equation}\label{coeff}
 \D_{\zbar}\Phi = \mu \D_z \Phi + \nu \overline{\D_z \Phi} \qquad \text{and} \qquad \D_{\zbar}\Psi = \mu \D_z \Psi + \nu \overline{\D_z \Psi},
\end{equation}
almost everywhere in $\Omega$. In the {\em regular set} $\cR_{\cF}$ of $\cF$\!, i.e., the set of points $z\in\Omega$ where the matrix
$$
M(z) = 
\left[ \begin{array}{cc}
\D_z \Phi(z) & \overline{\D_z \Phi(z)}  \\
\D_z \Psi(z) & \overline{\D_z \Psi(z)} 
\end{array} \right]
$$
is invertible, the values $\mu(z)$ and $\nu(z)$ are uniquely determined by \eqref{coeff}, that is, 
\begin{align}\label{coeff2}
 \mu(z) &= i\,\frac{\Psi_{\zbar}(z)\overline{\Phi_z(z)} - \overline{\Psi_z(z)}\Phi_{\zbar}(z)}{2\Im\bigl(\Phi_z(z) \overline{\Psi_z(z)}\bigr)}, \\\label{coeff3}
 \nu(z) &= i\,\frac{\Phi_{\zbar}(z) \Psi_z(z) - \Phi_z(z)\Psi_{\zbar}(z)}{2\Im\bigl(\Phi_z(z) \overline{\Psi_z(z)}\bigr)}.
\end{align}
Note that changing the generators corresponds to multiplying $M(z)$ by an invertible constant matrix. Hence the regular set and its complement, the {\em singular set}
$$
\cS_{\cF} = \bigl\{z \in \Omega :  2i \Im\bigl(\Phi_z(z) \overline{\Psi_z(z)}\,\bigr) = \det M(z) = 0\bigr\},
$$
depend only on the family $\cF$ and not the choice of generators.

It can be proven that for almost every $z\in \cS_{\cF}$ the vector $\bigl(\Phi_{\zbar}(z), \Psi_{\zbar}(z)\bigr)$ lies in the range of the linear operator $M(z) : \C^2 \to \C^2$. It follows that on the singular set one may define $\nu(z) = 0$. Here the assumption that the family $\cF$ consists entirely of quasiregular mappings is needed. By quasiregularity, one has for every $\alpha, \beta\in\R$
\begin{equation}\label{qrfamily}
|\alpha\, \D_{\zbar}\Phi(z) + \beta\, \D_{\zbar} \Psi(z)| \leq k|\alpha\, \D_z \Phi(z) + \beta\, \D_z \Psi(z)|, \qquad \text{for a.e. $z\in\Omega$}.
\end{equation}
There is a technical difficulty: the set where the inequality \eqref{qrfamily} holds depends, in general, on $\alpha$ and $\beta$. A short argument shows that \eqref{qrfamily} holds on the same set of full measure for all reals, see Lemma 12.1 in \cite{G-cpt} or p. 465 in \cite{AIM}.

Finally, ellipticity bounds in \eqref{geneq} follow for the singular set $\cS_{\cF}$ by definition of $\mu$ and $\nu$, since $\Phi$ and $\Psi$ are $K$-quasiregular. For the regular set one tests the inequality \eqref{qrfamily} by real-valued measurable functions $\theta(z)$ instead of parameters $\alpha$ and $\beta$. 

\medskip

As seen above, the existence of the general Beltrami equation \eqref{geneq} follows from local properties. For the uniqueness we need also global qualities. Here Theorem \ref{independent} comes into play. The coefficients are uniquely determined on the regular set $\cR_{\cF}$ by \eqref{coeff2} and \eqref{coeff3}. Moreover, Theorem \ref{independent} shows that the singular set $\cS_{\cF}$ has measure zero, thus proving the uniqueness.

\bigskip

\noindent {\em Step 2. General linear family}.\; It is a straightforward calculation after the two-dimensional case to achieve the same for general linear families.

Since the family consists of $K$-quasiregular mappings, we have the inequality
\begin{equation}\label{qrfamily2}
\left|\,\sum_{i\in\cI} a_i\, \partial_{\zbar}f_i(z)\right| \leq k\left|\,\sum_{i\in\cI} a_i\, \partial_z f_i(z)\right|, \qquad \text{for a.e. $z\in\Omega$},\quad k = \frac{K-1}{K+1}.
\end{equation}
There is the same technical difficulty as above; namely the set where the inequality \eqref{qrfamily2} holds depends, in general, on real numbers $a_i$. Noting that there are only countable many generators by assumption, the same argument as before works, see Lemma 12.1 in \cite{G-cpt} or p. 465 in \cite{AIM}. Hence we have a set $E \subset \Omega$ of full measure such that \eqref{qrfamily2} holds for all real coefficients $a_i$.

Note that two generating mappings, for example, $f_1$ and $f_2$, define a two-dimensional linear family of quasiregular mappings, and thus functions of the two-dimensional family satisfy the unique general linear Beltrami equation \eqref{gen}. Further, by Theorem \ref{independent}, we have a set of full measure $E'\subset E$ such that $\partial_z f_1(z)$ and $\partial_z f_2(z)$ are $\R$-linearly independent on $E'$.

\medskip

Our goal is to find $\cH:\Omega\times\C\to\C$ that satisfies
\begin{enumerate}
\item[(H1)] For $w_1, w_2\in\C$,
$$
|\cH(z,w_1)-\cH(z,w_2)|\leq k|w_1-w_2|, \qquad \text{for almost every $z\in\Omega$.}
$$
\item[(H2)] $\cH(z,0)\equiv0$.
\end{enumerate}
Moreover, we want every mapping $f\in\cF$ to solve the Beltrami equation
\begin{equation}\label{Hqr3}
\partial_{\zbar} f(z)=\cH(z,\partial_{z} f(z)).
\end{equation}

We define $\cH$ for $z\in E'$ by \eqref{Hqr3} going through all mappings $f\in\cF$. The function $\cH$ is not over-determined. Indeed, on $E'$ the inequality \eqref{qrfamily} holds, and hence if $\partial_z f(z) = \partial_z g(z)$ for some $f,g\in\cF$, then $\partial_{\zbar}f(z) = \partial_{\zbar}g(z)$. Moreover, since $\partial_z f_1(z)$ and $\partial_z f_2(z)$ are $\R$-linearly independent on $E'$, by our above remark, $\cH(z, w)$ is defined for all $w\in\C$.

The definition of $\cH$ and quasiregularity of mappings imply $k$-Lipschitz property on the second variable, that is, condition (H1). Also (H2) and that the equation \eqref{Hqr3} holds for all $f\in\cF$ follow straight from definition.

\medskip

Finally, we see that the linearity of the family $\cF$ is inherited by $\cH$, that is, $w \mapsto \cH(z, w)$ is $\R$-linear. Thus we have
$\cH(z,w) = \mu(z)w + \nu(z)\wbar$, where $|\mu(z)| + |\nu(z)|  \leq k < 1$. Since $\mu$ and $\nu$ are uniquely defined and measurable for two-dimensional linear families, our claim follows.
\end{proof}

\section*{Acknowledgements}

The author thanks Kari Astala and Tadeusz Iwaniec for stimulating discussions on the subject of this paper.

\medskip

\footnotesize{
\noindent Department of Mathematics and Statistics, 

\noindent P.O. Box 68, FI-00014 University of Helsinki, Finland

\smallskip

\noindent \texttt{jarmo.jaaskelainen@helsinki.fi}}

\end{document}